\documentclass[final,1p,times]{amsart}

\usepackage{amsmath,amsfonts,amssymb}

\usepackage{amsthm}
\usepackage[all]{xy}

\newcommand{\Ext}{\operatorname{Ext}}

\newcommand{\Hom}{\operatorname{Hom}}

\renewcommand{\mod}{\operatorname{mod}}
\newcommand{\umod}{\operatorname{\underline{mod}}}
\newcommand{\Triv}{\operatorname{T}}
\newcommand{\soc}{\operatorname{soc}}
\newcommand{\charact}{\operatorname{char}}

\newcommand{\Tr}{\operatorname{Tr}}
\newcommand{\op}{\operatorname{op}}

\newtheorem{theorem}{Theorem}[section]

\newtheorem{corollary}[theorem]{Corollary}
 \normalbaselines
\begin{document}

\title[Symmetric algebras of polynomial growth]
{Derived equivalence classification of symmetric algebras of
polynomial growth}

\author{Thorsten Holm}
\address{Institut f\"{u}r Algebra, Zahlentheorie und Diskrete
  Mathematik, Fa\-kul\-t\"{a}t f\"{u}r Ma\-the\-ma\-tik und Physik, Leibniz
  Universit\"{a}t Hannover, Welfengarten 1, 30167 Hannover, Germany}
%{Institut f\"{u}r Algebra und Geometrie, Fakult\"{a}t f\"{u}r
%Mathematik, Otto-von-Guericke-Universit\"{a}t Magdeburg, Postfach
%4120, 39016 Magdeburg, Germany\newline
%and\newline
%Department of Pure Mathematics, University of Leeds,
%Leeds LS2 9JT, United Kingdom}
\email{holm@math.uni-hannover.de}
\urladdr{http://www.iazd.uni-hannover.de/\~{ }tholm}

\author{Andrzej Skowro\'{n}ski}
\address{Faculty of Mathematics and Computer Science, 
    Nicolaus Copernicus University,
    Chopina 12/18, 87-100 Toru\'{n}, Poland}
\email{skowron@mat.uni.torun.pl}

%\author{Next author goes here}
%\address{Next author's postal address goes here}
%\email{Next author's mail address goes here}

\thanks{T.H. is supported by the research grant HO 1880/4-1
of the Deutsche Forschungsgemeinschaft (DFG), in the framework of 
the Research Priority Program SPP 1388 Representation Theory.}

\thanks{A.S. is supported by the research grant No. N N 201 269135 
of the Polish Ministry of Science and Higher Education.}



\begin{abstract} 
We complete the derived equivalence classification of all symmetric algebras
of polynomial growth, by solving the subtle problem of distinguishing the
standard and nonstandard nondomestic symmetric algebras of polynomial growth
up to derived equivalence. 
\end{abstract}

\maketitle

%%%%%%%%%%%%%%%%%%%%%%%%%%%%%%%%%%%%%%%%%%%%%%%%%%%%%%%%%%%%%%%%%%%%%%%%%%%%%%%
\section*{Introduction and the main result}
Throughout the article, $K$ will denote a fixed algebraically
closed field.
By an algebra is meant an associative, finite-dimensional $K$-algebra
with an identity.
For an algebra $A$, we denote by $\mod A$ the category of finite 
dimensional right $A$-modules and by $D$ the standard duality $\Hom_K(-,K)$
on $\mod A$.
An algebra $A$ is called \emph{selfinjective} if $A_A$ is an injective 
$A$-module, or equivalently, the projective $A$-modules are injective.
Prominent classes of selfinjective algebras are formed by the
\emph{Frobenius algebras} $A$ for which there exists an associative,
nondegenerate, $K$-bilinear form $(-,-) : A \times A \to K$, and the
\emph{symmetric algebras} $A$ for which there exists an associative,
symmetric, nondegenerate, $K$-bilinear form $(-,-) : A \times A \to K$.
By the classical theorems of T. Nakayama \cite{Nak1}, \cite{Nak2}, 
an algebra $A$ is Frobenius (respectively, symmetric) if and only if
$A \cong D(A)$ in $\mod A$ (respectively, as $A$-$A$-bimodules).
We also mention that every selfinjective algebra $A$ is Morita
equivalent to a Frobenius algebra, namely to its basic algebra.
Moreover, for every algebra $B$, the trivial extension
$\Triv(B) = B \ltimes D(B)$ of $B$ by the $B$-$B$-bimodule $D(B)$
is a symmetric algebra, and $B$ is a factor algebra of $\Triv(B)$.
It follows also from a result of T. Nakayama \cite{Nak2}
that the left socle and the right socle of a selfinjective algebra $A$
coincide, and we denote them by $\soc(A)$.
Two selfinjective algebras $A$ and $\Lambda$ are said to be
\emph{socle equivalent} if the factor algebras $A/\soc(A)$ 
and $\Lambda/\soc(\Lambda)$ are isomorphic.

According to the remarkable Tame and Wild Theorem of Y.A. Drozd \cite{Dr}
the class of (finite-dimensional) $K$-algebras over $K$ may be divided
into two disjoint classes.
One class consists of the tame algebras for which the indecomposable modules 
occur, in each dimension $d$, in a finite number of discrete and a finite
number of one-parameter families.
The second class consists of the wild algebras for which the representation
theory comprises the representation theories of all finite-dimensional algebras
over $K$ (see \cite[Chapter XIX]{SS2}).
Hence a classification of finite-dimensional modules is only feasible for tame
algebras.
More precisely, following Y.A. Drozd \cite{Dr}, an algebra $A$ is said to be
\emph{tame}, if for any positive integer $d$, there exists a finite number of
$K[x]$-$A$-bimodules $M_i$, $1 \leq i \leq n_d$, which are finitely generated 
and free as left $K[x]$-modules ($K[x]$ is the polynomial algebra in one variable 
over $K$) and all but finitely many isomorphism classes of indecomposable modules
of dimension $d$ in $\mod A$ are of the form $K[x]/(x-\lambda) \otimes_{K[x]} M_i$
for some $\lambda \in K$ and some $i \in \{1,\dots,n_d\}$.
Let $\mu_A(d)$ be the least number of $K[X]$-$A$-bimodules satisfying the above
condition for $d$.
Then $A$ is said to be of \emph{polynomial growth} (respectively, \emph{domestic})
if there exists a positive integer $m$ such that $\mu_A(d) \leq d^m$ (respectively, 
$\mu_A(d) \leq m $) for all $d \geq 1$. 
Moreover, from the validity of the second Brauer--Thrall conjecture, $\mu_A(d) = 0$ 
for all $d \geq 1$ if and only if $A$ is \emph{representation-finite} (there are 
only finitely many isomorphism classes of indecomposable modules in $\mod A$).

One central problem of modern representation theory is the determination 
of the module categories $\mod A$ of tame selfinjective algebras $A$.
Recently, the module categories of all selfinjective algebras of polynomial
growth have been described completely.
It has been proved by the second named author \cite{Sk3} that a nonsimple basic 
connected selfinjective algebra $A$ is of polynomial growth if and only if $A$
is socle equivalent to an orbit algebra $\widehat{B}/G$, where $\widehat{B}$
is the repetitive category of an algebra $B$, being a tilted algebra of
Dynkin or Euclidean type or a tubular algebra, and $G$ is an admissible infinite
cyclic automorphism group of $\widehat{B}$.
In particular, the Morita equivalence classification of the selfinjective algebras
of polynomial growth splits into two cases: the \emph{standard algebras} whose basic
algebras admit simply connected Galois coverings, and the remaining \emph{nonstandard
algebras} (see \cite{Sk1}).
We refer to the survey article \cite{Sk2} for the Morita equivalence classification
and the structure of module categories of the selfinjective algebras of polynomial
growth.

In this paper, we are concerned with the problem of derived equivalence classification
of selfinjective algebras of polynomial growth. For an algebra $A$, we denote by
$D^b(\mod A)$ the derived category of bounded complexes from $\mod A$.
Then two algebras $A$ and $\Lambda$ are said to be \emph{derived equivalent}
if the derived categories $D^b(\mod A)$ and $D^b(\mod \Lambda)$ are equivalent 
as triangulated categories.
Since Happel's work \cite{Ha1} interpreting tilting theory in terms of 
equivalences of derived categories, the machinery of derived categories
has been of interest to representation theorists.
In \cite{Ric1} J. Rickard proved his celebrated criterion: two algebras 
$A$ and $\Lambda$ are derived equivalent if and only if $\Lambda$ is 
the endomorphism algebra of a tilting complex over $A$. 
Since a lot of interesting properties are preserved by derived equivalences
(see Section \ref{sec:invariants}), it is for many purposes important to classify 
classes of algebras up to derived equivalence, instead of Morita equivalence.
In particular, for selfinjective algebras the representation types
introduced above are invariants of the derived category.

In \cite[Theorem 2.2]{Asa} H. Asashiba proved that the derived equivalence
classes of connected repre\-sen\-ta\-tion-finite standard (respectively, nonstandard)
selfinjective algebras are determined by the combinatorial data called the types,
and the derived equivalence classes of the standard and nonstandard 
representation-finite selfinjective algebras are disjoint.

A complete 
derived equivalence classification of the representation-infinite 
domestic standard (respectively, nonstandard) symmetric algebras has been 
established in our joint papers with R. Bocian \cite{BoHS1} (respectively, 
\cite{BoHS2}).
In \cite{HS} we completed the classification by showing that the derived
equivalence classes of the standard and nonstandard representation-infinite
domestic symmetric algebras are disjoint.
In fact, we established in \cite{HS}
the derived equivalence classification of all connected domestic
symmetric algebras by bound quiver algebras.

A 
%complete 
derived equivalence classification of the nondomestic standard 
(respectively, nonstandard) symmetric algebras of polynomial growth has been 
established in our joint papers with J. Bia{\l}kowski \cite{BiHS1} (respectively, 
\cite{BiHS2}). This classification is complete up to certain scalar parameters
occurring in the relations for which it seems intractable by current methods
to decide for which scalars the corresponding algebras are derived equivalent. 

The main open question in the derived equivalence classification of 
nondomestic symmetric algebras of polynomial growth
has been to distinguish the standard algebras from the 
nonstandard algebras up to derived equivalence.   
%We note that the problem of distinction of the derived equivalence classes of
%standard and 
%nonstandard symmetric algebras of polynomial growth 
This is subtle because the stable
Auslander-Reiten quivers of all these algebras consist only of stable tubes
(see \cite{NS}, \cite{Sk1}, \cite{Sk2}).  

In this paper we solve this problem by proving the following main result and hence 
complete (up to the scalars mentioned above)
the derived equivalence classification of the
symmetric algebras of polynomial growth.
\medskip

\noindent
{\bf Main Theorem.}
{\em Let $A$ be a standard selfinjective algebra and $\Lambda$ be a basic, connected,
nonstandard, nondomestic, symmetric algebra of polynomial growth.
Then $A$ and $\Lambda$ are not derived equivalent.
}
\medskip

In particular, it follows that the derived equivalence classes
of the standard and nonstandard nondomestic symmetric algebras
of polynomial growth are disjoint.
For a list of explicit representatives of the derived equivalence classes
we refer to Section \ref{Sec:normalform} below. 

The crucial tool for proving the Main Theorem are so called 
K\"ulshammer ideals defined by B. K\"ulshammer in 
\cite{Ku1}, \cite{Ku2}, \cite{Ku3}, \cite{Ku4} 
for symmetric algebras of positive characteristic, which have been
shown by A. Zimmermann \cite{Zi} to be invariants of derived equivalences.
These invariants are suitable for our purposes since the nonstandard
nondomestic selfinjective algebras of polynomial growth occur only in
characteristic $2$ and $3$, see \cite{BiS2}.
Moreover, in our proof of the Main Theorem we apply also a recent result 
of D. Al-Kadi \cite{AlK} describing the dimensions of the second 
Hochschild cohomology spaces of the preprojective algebra of Dynkin type
$\mathbb{D}_4$ and its unique proper deformation (in the sense of \cite{BES}).

In fact, the Main Theorem completes also the derived equivalence 
classification of tame symmetric algebras $A$ with periodic modules 
(all modules in $\mod A$ without projective direct summands are periodic
with respect to the action of the syzygy operator $\Omega_A$).
Recall that $\Omega_A$
assigns to a module $M$ in $\mod A$ the kernel of a minimal projective
cover $P_A(M) \to M$.
Since for the symmetric algebras $A$, the second syzygy 
$\Omega_A^2$ is the Auslander--Reiten
translation $\tau_A = D \Tr$ (as functors on the stable module category 
$\umod A$), the class of tame symmetric algebras $A$ with $\Omega_A$-periodic
module categories coincides with the class of tame symmetric algebras $A$
for which the stable Auslander--Reiten quiver $\Gamma_A^s$ consists only 
of stable tubes.
It has been proved recently by K. Erdmann and A. Skowro\'nski \cite{ES2} that 
a nonsimple, basic, connected, symmetric algebra is tame with periodic modules
if and only if $A$ is of one of the forms:
a representation-finite symmetric algebra, a nondomestic symmetric algebra
of polynomial growth, or an algebra of quaternion type (in the sense of \cite{Er}).
We refer to \cite[Chapter VII]{Er} and \cite[Section 5]{ES1} (see also 
\cite[Section 8]{Sk2}) for a Morita equivalence classification of algebras 
of quaternion type.
The derived equivalence classification of all algebras of quaternion type has 
been established by the first author in \cite[Section 5]{Ho} (see also 
\cite[Propositions 5.4 and 5.8]{ES1}).

For basic background on the representation theory we refer to the books
\cite{ASS}, 
%\cite{Er}, \cite{Ha2}, \cite{KeZ}, \cite{Rin}, 
\cite{SS1}, \cite{SS2} and to the survey articles \cite{Sk2} and \cite{Y}.

%%%%%%%%%%%%%%%%%%%%%%%%%%%%%%%%%%%%%%%%%%%%%%%%%%%%%%%%
\section{Invariants of derived equivalences of algebras}
\label{sec:invariants}

The aim of this section is to present some properties of algebras which
are invariant under derived equivalences.

The following results of J. Rickard \cite[Corollary~5.3]{Ric3}
(symmetric case) and S. Al-Nofayee \cite{AlN} (selfinjective case) establish
invariance of the classes of symmetric algebras and selfinjective 
algebras under derived equivalences.

\begin{theorem}
\label{th:Al-Nofayee} %
Let $A$ and $\Lambda$ be derived equivalent algebras.
Then the following equivalences hold.
\begin{enumerate}
\item[{(i)}] $A$ is symmetric if and only if $\Lambda$ is symmetric.
\item[{(ii)}] $A$ is selfinjective if and only if $\Lambda$ is selfinjective.
\end{enumerate}
\end{theorem}

Recall that two algebras $A$ and $\Lambda$ are said to be \emph{stably
equivalent} if the stable module categories $\umod A$ and $\umod \Lambda$
(modulo projectives) are equivalent.
The following result proved by J. Rickard in \cite[Corollary 2.2]{Ric2}
(see also \cite[Corollary 5.5]{Ric3})
is fundamental for the study of stable and derived equivalences
of selfinjective algebras. 

\begin{theorem}
\label{th:Ric} %
Let $A$ and $\Lambda$ be derived equivalent selfinjective algebras.
Then $A$ and $\Lambda$ are stably equivalent.
\end{theorem}

This together with the following theorem of H. Krause and G. Zwara 
\cite{KrZ} shows that the hierarchy of tame selfinjective
algebras is preserved by derived equivalences.

\begin{theorem}
\label{th:KrZ} %
Let $A$ and $\Lambda$ be stably equivalent algebras.
Then $A$ is tame (respectively, domestic, of polynomial growth)
if and only if $\Lambda$ has the same property.
\end{theorem}

For a selfinjective algebra $A$, denote by $\Gamma_A^s$ the
\emph{stable Auslander--Reiten quiver} of $A$, obtained from the
Auslander--Reiten quiver $\Gamma_A$ of $A$ by removing the projective
vertices and the arrows attached to them.

We have the following important consequence of Theorem \ref{th:Ric}.

\begin{corollary}
%\label{cor:Ric} %
\label{cor:1.4} %
Let $A$ and $\Lambda$ be derived equivalent selfinjective algebras.
Then $\Gamma_A^s$ and $\Gamma_{\Lambda}^s$ are isomorphic as translation 
quivers.
\end{corollary}

For an algebra $A$, consider the enveloping algebra 
$A^e = A^{\op} \otimes_K A$ of $A$.
Recall that the category $\mod A^e$ of finite dimensional right 
$A^e$-modules is equivalent to the category of finite dimensional 
$A$-$A$-bimodules.
Moreover, $A$ is selfinjective if and only if $A^e$ is selfinjective.
The algebra $A$ is a right $A^e$-module, via $a (x \otimes y) = x a y$
for $a \in A$, $x \in A^{\op}$, $y \in A$.
We may then consider the \emph{Hochschild cohomology algebra}
$$
    HH^*(A) = \Ext_{A^e}^*(A,A) = \bigoplus_{i \geq 0} \Ext_{A^e}^i(A,A)
$$
of $A$, which is a graded commutative $K$-algebra with respect to the
Yoneda product (see \cite[Corollary 1]{Gerstenhaber}).
%(see \cite{CE}, \cite{Ha3}, \cite{Hoc}).
We note that $HH^0(A)$ is isomorphic to the center $Z(A)$ of $A$.

The following theorem proved by J. Rickard in \cite[Proposition 2.5]{Ric3}
(see also \cite[Theorem 4.2]{Ha2} for a special case) 
shows that Hochschild cohomology yields invariants under derived 
equivalences.

\begin{theorem}
\label{th:1.5} %
Let $A$ and $B$ be derived equivalent algebras.
Then $HH^*(A)$ and $HH^*(B)$ are isomorphic as graded $K$-algebras.
\end{theorem}

\begin{corollary}
\label{cor:1.6} %
Let $A$ and $B$ be derived equivalent algebras.
Then the centers $Z(A)$ and $Z(B)$ are isomorphic as $K$-algebras.
\end{corollary}

We also note the following fact (see \cite[Lemma III.1.2]{Ha2}).

\begin{theorem}
\label{th:1.7} %
Let $A$ and $B$ be derived equivalent algebras.
Then the Grothendieck groups $K_0(A)$ and $K_0(B)$ are isomorphic.
\end{theorem}

%%%%%%%%%%%%%%%%%%%%%%%%%%%%%%%%%%%%%%%%%%%%%%%%%%%%%%%%%%%%%%%
\section{K\"ulshammer ideals}
\label{sec:Reynolds}

A prominent role in the proof of our main result is played by 
K\"ulshammer ideals, introduced by B. K\"ulshammer 
for symmetric algebras over algebraically closed fields of positive
characteristic. These form a decreasing sequence of ideals of the
center of an algebra, and this sequence 
has been shown by A. Zimmermann \cite{Zi} to be 
invariant under derived equivalences.
We recall the definition of K\"ulshammer ideals below; 
for more details on these invariants we refer to \cite{Ku1}, \cite{Ku2}, 
\cite{Ku3}, \cite{Ku4}, \cite{Ku5}, \cite{HHKM}, \cite{HZ}.

Let $K$ be an algebraically closed field of characteristic $p > 0$ 
and $A$ be a % finite dimensional 
symmetric $K$-algebra, i.e.
there exists an associative, symmetric, nondegenerate $K$-bilinear 
form $(-,-) : A \times A \to K$. 
For a $K$-subspace $M$ of $A$, denote by
$M^{\bot}$ the orthogonal complement of $M$ inside $A$ with
respect to the form $(-,-)$. Moreover, let $K(A)$ be the $K$-subspace
of $A$ generated by all commutators $[a,b]:=a b - b a$, for any
$a, b \in A$. For any $n \geq 0$ set
$$T_n (A) = \left\{ x \in A \mid x^{p^n} \in K(A)\right\}.$$
Then, by \cite{Ku1}, \cite{Ku2}, \cite{Ku3}, \cite{Ku4}, the orthogonal
complements $T_n (A)^{\bot}, n \geq 0$, are ideals of the center
$Z(A)$ of $A$, called \emph{K\"ulshammer ideals}.
They form a descending sequence
$$Z(A) = T_0(A)^{\perp} \supseteq T_1(A)^{\perp} \supseteq T_2(A)^{\perp}
\supseteq T_3(A)^{\perp} \supseteq \ldots$$
In fact, B.\,K\"ulshammer proved in \cite{Ku5} that the
equation $(\xi_n(z),x)^{p^n} = (z,x^{p^n})$ for any $x, z \in
Z(A)$ defines a mapping $\xi_n : Z(A) \rightarrow Z(A)$ whose image
$\xi_n(Z(A)) = T_n (A)^{\bot}$ is precisely the $n$-th K\"ulshammer
ideal.

Then we have the following theorem proved recently by
A.\,Zimmermann \cite[Theorem~1]{Zi} showing that the
K\"ulshammer ideals are derived invariants.

\begin{theorem}
\label{th:zimmermann} %
Let $A$ and $B$ be derived equivalent symmetric algebras over an
algebraically closed field of positive characteristic $p$. Then
there is an isomorphism $\varphi : Z(A) \rightarrow Z(B)$ of the
centers of $A$ and $B$ such that $\varphi(T_n (A)^{\bot}) = T_n
(B)^{\bot}$ for all nonnegative integers $n$.
\end{theorem}

Hence the sequence of K\"ulshammer ideals gives new
derived invariants, for symmetric algebras over algebraically
closed fields of positive characteristic.

%\smallskip

%In the next section we shall use this invariant for proving our main result. 

For using K\"ulshammer ideals in the context of derived equivalence 
classifications one has to be able to perform explicit computations 
with them and for this we shall later give 
explicit symmetrizing bilinear forms on the symmetric algebras under
consideration. 

\section{Derived normal forms of nondomestic symmetric algebras 
    of polynomial growth} \label{Sec:normalform}

In this section we present derived normal forms, i.e.
representatives of the derived equivalence classes, of the connected
nondomestic symmetric algebras of polynomial growth.

\medskip

Consider the following families of bound quiver algebras
(in the notation of \cite{BiHS2}):
\medskip

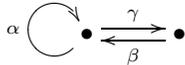
\begin{minipage}[t]{6cm}
$\Lambda_2$:
$$
  \xymatrix{
    \bullet \ar `ld_u[] `_rd[]^{\alpha} [] 
      \ar@<.5ex>^{\gamma}[r]
      & \bullet \ar@<.5ex>^{\beta}[l] &
  }
$$
\center \vspace{-.5pc}
 $\alpha^2 \gamma = 0$, $\beta \alpha^2 = 0$,
 $\gamma \beta \gamma = 0$, $\beta \gamma \beta = 0$,
 $\beta \gamma = \beta \alpha \gamma$,
 $\alpha^3 = \gamma \beta$
\end{minipage} \hfill
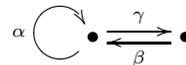
\begin{minipage}[t]{6cm}
$\Lambda_2'$:
$$
  \xymatrix{
    \bullet \ar `ld_u[] `_rd[]^{\alpha} [] 
      \ar@<.5ex>^{\gamma}[r]
      & \bullet \ar@<.5ex>^{\beta}[l] &
  }
$$
\center \vspace{-.5pc} 
 $\alpha^3 = \gamma \beta$,
 $\beta \gamma = 0$, \\
 $\beta \alpha^2 = 0$,
 $\alpha^2 \gamma = 0$
\end{minipage}

\bigskip

\begin{minipage}[t]{6cm}
$\Lambda_3(\lambda)$, $\lambda \in K \setminus \{0,1\} $:
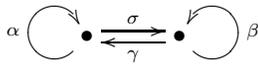
$$
  \xymatrix{
%    \bullet \ar `ru_d[] `_lu[]^{\alpha} [] 
    \bullet \ar `ld_u[] `_rd[]^{\alpha} [] 
      \ar@<.5ex>^{\sigma}[r]
      &     \bullet \ar `rd^u[] `^ld[]_{\beta} [] 
           \ar@<.5ex>^{\gamma}[l] &
  }
$$
\center \vspace{-.5pc} 
    $\alpha^4 = 0$, $\gamma \alpha^2 = 0$,
    $\alpha^2 \sigma = 0$,
    $\alpha^2 = \sigma \gamma + \alpha^3$,
    $\lambda \beta^2 = \gamma \sigma$,
    $\gamma \alpha = \beta \gamma$,
    $\sigma \beta = \alpha \sigma$
\end{minipage} \hfill
\begin{minipage}[t]{6cm}
$\Lambda'_3(\lambda)$, $\lambda \in K \setminus \{0,1\} $:
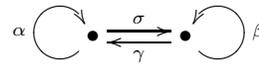
$$
  \xymatrix{
%    \bullet \ar `ru_d[] `_lu[]^{\alpha} [] 
    \bullet \ar `ld_u[] `_rd[]^{\alpha} [] 
      \ar@<.5ex>^{\sigma}[r]
      &     \bullet \ar `rd^u[] `^ld[]_{\beta} [] 
           \ar@<.5ex>^{\gamma}[l] &
  }
$$
\center \vspace{-.5pc} 
$\alpha^2 = \sigma \gamma$,
$\lambda \beta^2 = \gamma \sigma$,\\
$\gamma \alpha = \beta \gamma$,
$\sigma \beta = \alpha \sigma$
\end{minipage}

\bigskip

\begin{minipage}[t]{6cm}
$\Lambda_5$:
$$
  \xymatrix{
%    \bullet \ar `ru_d[] `_lu[]^{\alpha} [] 
    \bullet 
      \ar@<.5ex>^{\beta}[r]
      & \bullet \ar@(lu,ur)[]^{\alpha}
          \ar@<.5ex>^{\gamma}[l] 
      \ar@<.5ex>^{\delta}[r]
      & \bullet \ar@<.5ex>^{\sigma}[l] 
  }
$$
\center \vspace{-.5pc}
$\alpha^2 = \gamma \beta$,
$\alpha^3 = \delta \sigma$,
$\beta \delta = 0$,
$\sigma \gamma = 0$,
$\alpha \delta = 0$,
$\sigma \alpha = 0$,
$\gamma \beta \gamma = 0$,
$\beta \gamma \beta = 0$,
$\beta \gamma = \beta \alpha \gamma$
\end{minipage} \hfill
\begin{minipage}[t]{6cm}
$\Lambda_5'$:
$$
  \xymatrix{
%    \bullet \ar `ru_d[] `_lu[]^{\alpha} [] 
    \bullet 
      \ar@<.5ex>^{\beta}[r]
      & \bullet \ar@(lu,ur)[]^{\alpha}
          \ar@<.5ex>^{\gamma}[l] 
      \ar@<.5ex>^{\delta}[r]
      & \bullet \ar@<.5ex>^{\sigma}[l] 
  }
$$
\center \vspace{-.5pc}
$\alpha^2 = \gamma \beta$,
$\beta \delta = 0$,
$\beta \gamma = 0$,
$\sigma \gamma = 0$,
$\alpha \delta = 0$,
$\sigma \alpha = 0$,
$\alpha^3 = \delta \sigma$
\end{minipage} 

\bigskip

\begin{minipage}[t]{6cm}
$\Lambda_9$:
$$
  \xymatrix@R=.9pc{
    & \bullet \ar@<.5ex>^{\gamma}[dd] \\ \\
    & \bullet \ar@<.5ex>^{\delta}[uu]
      \ar@<.5ex>^{\beta}[ld]
      \ar@<.5ex>^{\varepsilon}[rd] \\
    \bullet \ar@<.5ex>^{\alpha}[ur]
      && \bullet \ar@<.5ex>^{\xi}[ul] \\
  }
$$
\center \vspace{-.5pc}
$\beta \alpha + \delta \gamma + \varepsilon \xi = 0$,
$\gamma \delta = 0$,
$\xi \varepsilon = 0$,
$\alpha \beta \alpha = 0$,
$\beta \alpha \beta = 0$,
$\alpha \beta = \alpha \delta \gamma \beta$
\end{minipage} \hfill
\begin{minipage}[t]{6cm}
$\Lambda_9'$:
$$
  \xymatrix@R=.9pc{
    & \bullet \ar@<.5ex>^{\gamma}[dd] \\ \\
    & \bullet \ar@<.5ex>^{\delta}[uu]
      \ar@<.5ex>^{\beta}[ld]
      \ar@<.5ex>^{\varepsilon}[rd] \\
    \bullet \ar@<.5ex>^{\alpha}[ur]
      && \bullet \ar@<.5ex>^{\xi}[ul] \\
  }
$$
\center \vspace{-.5pc}
$\beta \alpha + \delta \gamma + \varepsilon \xi = 0$,
$\alpha \beta = 0$,
$\xi \varepsilon = 0$,
$\gamma \delta = 0$
\end{minipage}

\bigskip

\begin{minipage}[t]{6cm}
$A_1(\lambda)$, $\lambda \in K \setminus \{0,1\}$:
$$
  \xymatrix{
    \bullet 
      \ar@<.5ex>^{\alpha}[r]
      & \bullet \ar@<.5ex>^{\gamma}[l] 
      \ar@<.5ex>^{\sigma}[r]
      & \bullet \ar@<.5ex>^{\beta}[l] 
  }
$$
\center \vspace{-.5pc}
$\alpha \gamma \alpha = \alpha \sigma \beta$,
$\beta \gamma \alpha = \lambda \beta \sigma \beta$,
$\gamma \alpha \gamma = \sigma \beta \gamma$, 
$\gamma \alpha \sigma = \lambda \sigma \beta \sigma$
\end{minipage} \hfill
\begin{minipage}[t]{6cm}
$A_4$:
$$
  \xymatrix@R=.9pc{
    & \bullet \ar@<.5ex>^{\gamma}[dd] \\ \\
    & \bullet \ar@<.5ex>^{\delta}[uu]
      \ar@<.5ex>^{\beta}[ld]
      \ar@<.5ex>^{\varepsilon}[rd] \\
    \bullet \ar@<.5ex>^{\alpha}[ur]
      && \bullet \ar@<.5ex>^{\xi}[ul] \\
  }
$$
\center \vspace{-.5pc}
$\beta \alpha + \delta \gamma + \varepsilon \xi = 0$,
$\alpha \beta = 0$,
$\gamma \varepsilon = 0$,
$\xi \delta = 0$
\end{minipage}

\bigskip

%\bigskip
%
%{\Large To be completed}
%
%\bigskip

Note that the algebra $\Lambda'_9$ is just the preprojective algebra
of Dynkin type $\mathbb{D}_4$.
\bigskip

It has been proved in \cite[Theorems 1 and 2]{BiS1} and 
\cite[Theorem 1.1]{BiS2} that:
\begin{itemize}
 \item
  $\Lambda_2$ and $\Lambda_2'$ are symmetric algebras of tubular type
  $(3,3,3)$,
  and $\Lambda_2 \cong \Lambda_2'$ if and only if $\charact K \neq 3$.
  Moreover, $\Lambda_2$ and $\Lambda_2'$ are socle equivalent.
 \item
  $\Lambda_3(\lambda)$ and $\Lambda_3'(\lambda)$, 
  $\lambda \in K \setminus \{0,1\}$,
  are symmetric algebras of tubular type $(2,2,2,2)$
  and $\Lambda_3(\lambda) \cong \Lambda_3'(\lambda)$ 
  if and only if $\charact K \neq 2$.
  Moreover, $\Lambda_3(\lambda)$ and $\Lambda_3'(\lambda)$ are socle equivalent.
 \item
  $\Lambda_5$ and $\Lambda_5'$ are symmetric algebras of tubular type
  $(2,4,4)$,
  and $\Lambda_5 \cong \Lambda_5'$ if and only if $\charact K \neq 2$.
  Moreover, $\Lambda_5$ and $\Lambda_5'$ are socle equivalent.
 \item
  $\Lambda_9$ and $\Lambda_9'$ are weakly symmetric algebras of tubular type
  $(3,3,3)$,
  and $\Lambda_9 \cong \Lambda_9'$ if and only if $\charact K \neq 2$.
  Furthermore, $\Lambda_9$ and $\Lambda_9'$ are symmetric algebras 
  if and only if $\charact K = 2$.
  Moreover, $\Lambda_9$ and $\Lambda_9'$ are socle equivalent.
 \item
  $A_1(\lambda)$, $\lambda \in K \setminus \{0,1\}$,
  are symmetric algebras of tubular type $(2,2,2,2)$.
 \item
  $A_4$ is a symmetric algebra of tubular type $(3,3,3)$.
\end{itemize}

It follows also from \cite[Theorem 2.1]{BiHS2} that $\Lambda_i'$
is a geometric degeneration of $\Lambda_i$ for $i \in \{2,5,9\}$, 
and $\Lambda_3'(\lambda)$ is a geometric degeneration 
of $\Lambda_3(\lambda)$, for $\lambda \in K \setminus \{0,1\}$.

Consider also the following family of the trivial extensions 
of Ringel's tubular canonical algebras:
%\begin{itemize}
% \item
%    $\Lambda_1(c)$, $c \in K \setminus \{ 0, 1 \}$, of tubular type
%    $(2,2,2,2)$ given by the bound quiver
%%\xy
%%\endxy

\medskip

$\Lambda(2,2,2,2,\lambda), \lambda \in K \setminus \{0,1\}$:
%\xyoption{all}
\newdir{m}{{}*!/-5pt/@{}}
\newdir{>m}{@{>}*!/-5pt/{}}
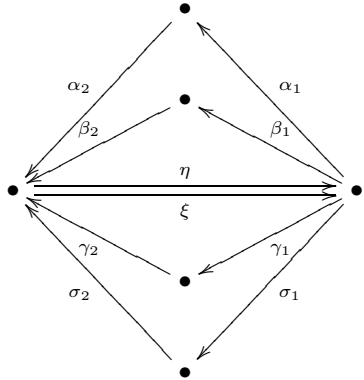
\[
\xymatrix{
 & & {\bullet} \ar[ddll]_{\alpha_2} & & \\
 & & {\bullet} \ar[dll]_{\beta_2} & & \\
 {\bullet} \ar@{m->m}@<.4ex>[rrrr]^{\eta} \ar@{m->m}@<-.4ex>[rrrr]_{\xi}
 & & & & {\bullet}  \ar[uull]_{\alpha_1} \ar[ull]_{\beta_1}
                    \ar[dll]^{\gamma_1} \ar[ddll]^{\sigma_1} \\
 & & {\bullet} \ar[ull]^{\gamma_2} & & \\
 & & {\bullet} \ar[uull]^{\sigma_2} & & \\
}
\]
    $\alpha_1 \alpha_2 +   \beta_1 \beta_2 + \gamma_1 \gamma_2 = 0$,
    $\alpha_1 \alpha_2 + \lambda \beta_1 \beta_2 + \sigma_1 \sigma_2 = 0$,
    $\eta \alpha_1 = 0$, $\alpha_2 \eta = 0$,
    $\xi \beta_1 = 0$, $\beta_2 \xi = 0$, $\eta \gamma_1 = \xi \gamma_1$,
    $\gamma_2 \eta = \gamma_2 \xi$,
    $\eta \sigma_1 = \lambda \xi \sigma_1$, $\sigma_2 \eta = \lambda \sigma_2 \xi$.

\medskip

% \item
%    $\Lambda_2$ of tubular type $(3,3,3)$ given by the bound quiver
$\Lambda(3,3,3)$:
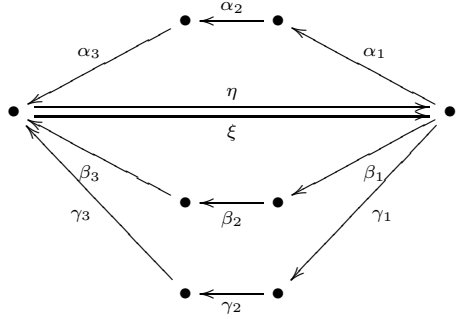
\[
\xymatrix{
 & & {\bullet} \ar[dll]_{\alpha_3} & {\bullet} \ar[l]_{\alpha_2} & & \\
 {\bullet} \ar@{m->m}@<.4ex>[rrrrr]^{\eta} \ar@{m->m}@<-.4ex>[rrrrr]_{\xi}
  & & & & & {\bullet}  \ar[ull]_{\alpha_1} \ar[dll]^{\beta_1}
                       \ar[ddll]^{\gamma_1} \\
 & & {\bullet} \ar[ull]^{\beta_3} & {\bullet} \ar[l]^{\beta_2} & & \\
 & & {\bullet} \ar[uull]^{\gamma_3} & {\bullet} \ar[l]^{\gamma_2} & & \\
}
\]
    $\alpha_1 \alpha_2 \alpha_3 + \beta_1 \beta_2 \beta_3 +
     \gamma_1 \gamma_2 \gamma_3 = 0$,
    $\eta \alpha_1 = 0$, $\alpha_3 \eta = 0$, $\xi \beta_1 = 0$,
    $\beta_3 \xi = 0$,
    $\eta \gamma_1 = \xi \gamma_1$, $\gamma_3 \eta = \gamma_3 \xi$,
    $\alpha_2 \alpha_3 \xi \alpha_1 \alpha_2 = 0$,
    $\beta_2 \beta_3 \eta \beta_1 \beta_2 = 0$,
    $\gamma_2 \gamma_3 \eta \gamma_1 \gamma_2 = 0$.

\medskip

% \item
%    $\Lambda_3$ of tubular type $(2,4,4)$ given by the bound quiver
$\Lambda(2,4,4)$:
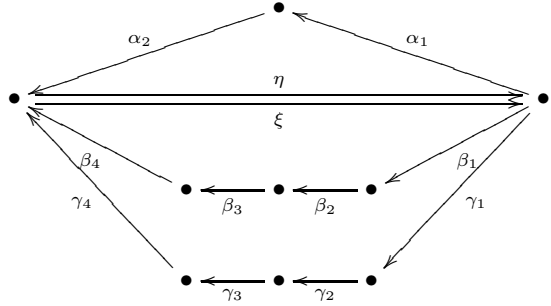
\[
\xymatrix{
 & & & {\bullet} \ar[dlll]_{\alpha_2} & & & \\
 {\bullet} \ar@{m->m}@<.35ex>[rrrrrr]^{\eta} \ar@{m->m}@<-.4ex>[rrrrrr]_{\xi}
  & & & & & & {\bullet}  \ar[ulll]_{\alpha_1} \ar[dll]^{\beta_1}
                       \ar[ddll]^{\gamma_1} \\
 & & {\bullet} \ar[ull]^{\beta_4} & {\bullet} \ar[l]^{\beta_3} &
  {\bullet} \ar[l]^{\beta_2} & & \\
 & & {\bullet} \ar[uull]^{\gamma_4} & {\bullet} \ar[l]^{\gamma_3} &
  {\bullet} \ar[l]^{\gamma_2} & & \\
}
\]
    $\alpha_1 \alpha_2 + \beta_1 \beta_2 \beta_3 \beta_4 +
     \gamma_1 \gamma_2 \gamma_3 \gamma_4 = 0$,  
    $\eta \alpha_1 = 0$, $\alpha_2 \eta = 0$, $\xi \beta_1 = 0$,
    $\beta_4 \xi = 0$,
    $\eta \gamma_1 = \xi \gamma_1$, $\gamma_4 \eta = \gamma_4 \xi$,
    $\beta_2 \beta_3 \beta_4 \eta \beta_1 \beta_2 = 0$,
    $\beta_3 \beta_4 \eta \beta_1 \beta_2 \beta_3 = 0$,
    $\gamma_2 \gamma_3 \gamma_4 \eta \gamma_1 \gamma_2 = 0$,
    $\gamma_3 \gamma_4 \eta \gamma_1 \gamma_2 \gamma_3 = 0$.

\medskip

% \item
%    $\Lambda_4$ of tubular type $(2,3,6)$ given by the bound quiver
$\Lambda(2,3,6)$:
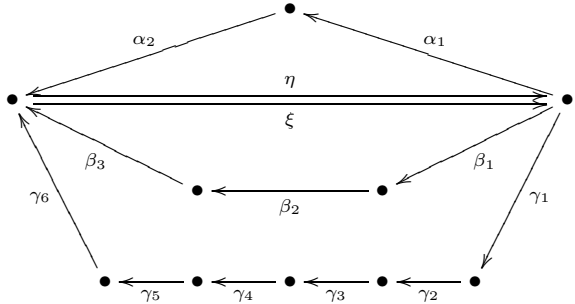
\[
\xymatrix{
 & & & {\bullet} \ar[dlll]_{\alpha_2} & & & \\
 {\bullet} \ar@{m->m}@<.35ex>[rrrrrr]^{\eta} \ar@{m->m}@<-.4ex>[rrrrrr]_{\xi}
  & & & & & & {\bullet}  \ar[ulll]_{\alpha_1} \ar[dll]^{\beta_1}
                       \ar[ddl]^{\gamma_1} \\
 & & {\bullet} \ar[ull]^{\beta_3} & & {\bullet} \ar[ll]^{\beta_2} & & \\
 & {\bullet} \ar[uul]^{\gamma_6} & {\bullet} \ar[l]^{\gamma_5} &
  {\bullet} \ar[l]^{\gamma_4} & {\bullet} \ar[l]^{\gamma_3} &
  {\bullet} \ar[l]^{\gamma_2} & \\
}
\]
    $\alpha_1 \alpha_2 + \beta_1 \beta_2 \beta_3 +
     \gamma_1 \gamma_2 \gamma_3 \gamma_4 \gamma_5 \gamma_6 = 0$,    
    $\eta \alpha_1 = 0$, $\alpha_2 \eta = 0$, $\xi \beta_1 = 0$,
    $\beta_3\xi = 0$,
    $\eta \gamma_1 = \xi \gamma_1$, $\gamma_6 \eta = \gamma_6 \xi$,
    $\beta_2 \beta_3 \eta \beta_1 \beta_2 = 0$,
     $\gamma_2 \gamma_3 \gamma_4 \gamma_5 \gamma_6 \eta \gamma_1 \gamma_2 = 0$,
     $\gamma_3 \gamma_4 \gamma_5 \gamma_6 \eta \gamma_1 \gamma_2 \gamma_3 = 0$,
     $\gamma_4 \gamma_5 \gamma_6 \eta \gamma_1 \gamma_2 \gamma_3 \gamma_4 = 0$,
     $\gamma_5 \gamma_6 \eta \gamma_1 \gamma_2 \gamma_3 \gamma_4 \gamma_5 = 0$.

%\end{itemize}

\medskip

A direct calculation shows that 
$\Lambda(2,2,2,2,\lambda) \cong \Triv(C(2,2,2,2,\lambda))$,
$\Lambda(3,3,3) \cong \Triv(C(3,3,3))$,
$\Lambda(2,4,4) \cong \Triv(C(2,4,4))$ and
$\Lambda(2,3,6) \cong \Triv(C(2,3,6))$,
where $C(2,2,2,2,\lambda)$, $C(3,3,3)$, $C(2,4,4)$ and $C(2,3,6)$
are the canonical algebras of tubular types
$(2,2,2,2)$, $(3,3,3)$, $(2,4,4)$ and $(2,3,6)$, respectively.

The following derived equivalence classification of the standard 
nondomestic symmetric algebras of polynomial growth follows from
\cite[Theorem]{BiHS1} (see also \cite{HR} for the derived equivalence 
of trivial extensions of tubular algebras) and the facts described
above.

\begin{theorem}
\label{th:3.1}
Let $A$ be a connected standard nondomestic symmetric algebra
of polynomial growth.
Then $A$ is derived equivalent to one of the following algebras:
\begin{itemize}
 \item
  two simple modules: $\Lambda_2'$ and $\Lambda_3'(\lambda), \lambda \in K \setminus \{0,1\}$;
 \item
  three simple modules: $\Lambda_5'$ and $A_1(\lambda), \lambda \in K \setminus \{0,1\}$;
 \item
  four simple modules: $\Lambda_9' \,(\charact K = 2)$ and $A_4$;
 \item
  six simple modules: $\Lambda(2,2,2,2,\lambda), \lambda \in K \setminus \{0,1\}$;
 \item
  eight simple modules: $\Lambda(3,3,3)$;
 \item
  nine simple modules: $\Lambda(2,4,4)$;
 \item
  ten simple modules: $\Lambda(2,3,6)$.
\end{itemize}
\end{theorem}

Note that in the cases of four, eight, nine and ten simple modules,
the above classification does not depend on a scalar.
Moreover, for any $\lambda \in K \setminus \{0,1\}$, the centers 
of the algebras $\Lambda_2'$ and $\Lambda_3'(\lambda)$ are not 
isomorphic (\cite[Lemma 2.2]{BiHS1}), and hence 
$\Lambda_2'$ and $\Lambda_3'(\lambda)$ are not derived equivalent,
by Corollary \ref{cor:1.6}. 
Similarly, the centers 
of the algebras $\Lambda_9'$ and $A_4$ are not 
isomorphic (\cite[Lemma 4.6]{BiHS1}), and hence 
$\Lambda_9'$ and $A_4$ are not derived equivalent. 
Moreover, for any $\lambda \in K \setminus \{0,1\}$,
the algebras $\Lambda_5'$ and $A_1(\lambda)$ have nonisomorphic
stable Auslander--Reiten quivers (see the proof of the Main Theorem),
and hence are not derived equivalent, by Corollary \ref{cor:1.4}.
Therefore, the classification is complete up to the scalars
$\lambda \in K \setminus \{0,1\}$ occuring in the algebras 
$\Lambda_3'(\lambda)$ and $A_1(\lambda)$.
We do not know how to decide for which scalars $\lambda, \mu \in K \setminus \{0,1\}$
the algebras $\Lambda_3'(\lambda)$ and $\Lambda_3'(\mu)$
(respectively, $A_1(\lambda)$ and $A_1(\mu)$)
are derived equivalent.

\medskip

It has been proved in \cite[Theorem 1.1]{BiS2} that the nonstandard
nondomestic symmetric algebras of polynomial growth occur only in
characteristic $2$ and $3$.
Furthermore, the following theorem proved in \cite{BiHS2}
gives the derived equivalence classification of these algebras
(when comparing with the results of \cite{BiHS2} please note that the algebra
$\Lambda_{10}$ occurring there is not symmetric and therefore does not have
to be considered here).

\begin{theorem}
\label{th:3.2}
Let $A$ be a connected nonstandard nondomestic symmetric algebra
of polynomial growth.
Then $A$ is derived equivalent to one of the following algebras:
\begin{itemize}
 \item
  two simple modules: $\Lambda_2 \, (\charact K = 3)$ and 
                      $\Lambda_3(\lambda), \lambda \in K \setminus \{0,1\} \, (\charact K = 2)$;
 \item
  three simple modules: $\Lambda_5 \, (\charact K = 2)$;
 \item
  four simple modules: $\Lambda_9 \, (\charact K = 2)$.
\end{itemize}
\end{theorem}

The above classification is complete up to the scalars 
$\lambda \in K \setminus \{0,1\}$ occuring in the algebras $\Lambda_3(\lambda)$.
We do not know how to decide for which scalars $\lambda, \mu \in K \setminus \{0,1\}$
and $K$ of characteristic $2$, the algebras $\Lambda_3(\lambda)$ and $\Lambda_3(\mu)$
are derived equivalent.

%%%%%%%%%%%%%%%%%%%%%%%%%%%%%%%%%%%%%%%%%%%%%%%%%%%%%%%%%%%%%%%%%%%%%
\section{Proof of the Main Theorem}
\label{sec:proof}

Let $\Lambda$ be a basic, connected, nonstandard, nondomestic, 
symmetric algebra of polynomial growth over $K$.
Then it follows from Theorem \ref{th:3.2} that either
$\charact K = 3$ and $\Lambda$ is derived equivalent to 
the algebra $\Lambda_2$ or $\charact K = 2$ and $\Lambda$ 
is derived equivalent to one of the algebras 
$\Lambda_3(\lambda), \, \lambda \in K \setminus \{0,1\}$, 
$\Lambda_5$, or $\Lambda_9$.

Assume now that $A$ is a standard selfinjective algebra over $K$
which is derived equivalent to $\Lambda$.
Then, by Theorems \ref{th:Al-Nofayee}, \ref{th:Ric} and \ref{th:KrZ},
$A$ is a nondomestic symmetric algebra of polynomial growth.
Moreover, by Corollary \ref{cor:1.6}, the centers $Z(A)$ and 
$Z(\Lambda)$ are isomorphic, and hence $A$ is connected, because $\Lambda$
is connected (by assumption on $\Lambda$).
Since $A$ is by assumption a standard algebra, we conclude that $A$
is derived equivalent to one of the algebras listed in Theorem \ref{th:3.1}.
Furthermore, it follows from Theorem \ref{th:1.7} 
that the Grothendieck groups $K_0(\Lambda)$ and $K_0(A)$ are isomorphic,
and hence the quivers $Q_{\Lambda}$ of $\Lambda$ and $Q_A$ of $A$ 
have the same number of vertices (see \cite[II.3 and III.3]{ASS}).
Finally, by the general theory of selfinjective algebras of tubular type
(see \cite[Section 5]{Sk2}, or \cite[Section 3]{Sk1}) we know that the
stable Auslander-Reiten quivers $\Gamma_{\Lambda_5}$ and $\Gamma_{\Lambda_5'}^s$
consist of $\mathbb{P}_1(K)$-families of stable tubes of tubular type $(2,4,4)$,
while the stable Auslander-Reiten quivers 
$\Gamma_{A_1(\lambda)}^s, \lambda \in K \setminus \{0,1\}$,
consist of $\mathbb{P}_1(K)$-families of stable tubes of tubular type $(2,2,2,2)$.
Hence, by Corollary \ref{cor:1.4}, $\Lambda_5$ (respectively, $\Lambda_5'$)
is not derived equivalent to an algebra 
$A_1(\lambda), \lambda \in K \setminus \{0,1\}$.

Moreover, $\Lambda_2$ is not derived equivalent to $\Lambda'_3(\lambda)$,
since the former is of tubular type $(3,3,3)$ and the latter is 
of tubular type $(2,2,2,2)$. 
Similarly, $\Lambda_3(\lambda)$ is not derived equivalent to $\Lambda'_2$. 

Finally, $\Lambda_9$ is not derived equivalent to $A_4$ since their
centers have different dimensions. In fact, the center of $\Lambda_9$
is of dimension 5, generated by the unit element $1$ and the (one-dimensional)
socles of
the four projective indecomposable modules; on the other hand, the center
of $A_4$ has an additional basis element $\beta\alpha-\xi\epsilon-\gamma\delta$
and hence has dimension 6 (cf. also \cite[Lemma 4.6]{BiHS1}).

\smallskip

Therefore, in order to prove the Main Theorem, it remains to show that:
\begin{enumerate}
 \item[{(1)}]
  For $\charact K = 3$, the algebras $\Lambda_2$ and $\Lambda_2'$
  are not derived equivalent;
 \item[{(2)}]
  For $\charact K = 2$ and $\lambda, \mu \in K \setminus \{0,1\}$,
  the algebras $\Lambda_3(\lambda)$ and $\Lambda_3(\mu)'$
  are not derived equivalent;
 \item[{(3)}]
  For $\charact K = 2$, the algebras $\Lambda_5$ and $\Lambda_5'$
  are not derived equivalent;
 \item[{(4)}] 
  For $\charact K = 2$, the algebras $\Lambda_9$ and $\Lambda_9'$
  are not derived equivalent.
\end{enumerate}

In order to prove (1)--(3) we will use the sequences of K\"ulshammer ideals, 
as described in Section \ref{sec:Reynolds}. 
To this end we have to provide below explicit associative,
nondegenerate symmetric bilinear forms on the algebras occurring
in (1)--(3).
It is well-known that such a form exists (i.e. that an algebra $A$
is symmetric) if and only if there is a $K$-linear form $\Psi:A\to K$,
called a symmetrizing form, 
such that $\Psi(ab)=\Psi(ba)$ for all $a,b\in A$ and the kernel
of $\Psi$ does not contain any nonzero left or right
ideal of $A$. 

The algebras
in (4) can not be distinguished by the K\"ulshammer ideals; for proving
(4) we shall use a recent result of D. Al-Kadi \cite{AlK} on the
second Hochschild cohomology of the two algebras involved.  
\smallskip
 
\begin{enumerate}  
\item[{(1)}]
  Let $\charact K = 3$. For both of the algebras $\Lambda_2$ and $\Lambda_2'$,
  a $K$-linear basis
 consisting of the different nonzero
paths is given by 
$$\mathcal{B} = \{ e_1,\alpha,\alpha^2,\alpha^3=\gamma\beta,\alpha^4=
\alpha\gamma\beta=\gamma\beta\alpha,\gamma,\alpha\gamma,e_2,\beta,
\beta\alpha,\beta\alpha\gamma \}.$$
Also for both algebras $\Lambda=\Lambda_2$ and $\Lambda=\Lambda'_2$, 
the center is as $K$-vector space generated as 
follows
$$Z(\Lambda) = \langle 1,\alpha^2,\alpha^3,\alpha^4,\beta\alpha\gamma
\rangle_K.
$$
It follows from \cite[Section 4]{BiS1} (or can be checked directly)
that a symmetrizing form $\Psi'$ on the standard algebra  
$\Lambda_2'$ is given by assigning 1 to $\alpha^4$ and $\beta\alpha\gamma$, 
and 0 to all other elements of the basis $\mathcal{B}$.

Note that we can not use the same form for the nonstandard algebra
$\Lambda_2$, since $\Psi'(\gamma\beta) = 0$ whereas 
$\Psi'(\beta\gamma)=\Psi'(\beta\alpha\gamma)=1$, i.e. the form
is not symmetric. Instead, a symmetrizing form $\Psi$
on $\Lambda_2$ 
is given by assigning 1 to $\alpha^4$, to $\beta\alpha\gamma$, and to 
$\alpha^3=\gamma\beta$, and  
assigning 0 to all other elements of the basis $\mathcal{B}$.

Now we have to examine the commutator subspaces. 

First, we consider the standard algebra $\Lambda_2'$. Since $\beta\gamma=0$ 
we obtain
$$K(\Lambda_2')=\langle \gamma,\beta,\alpha\gamma,\beta\alpha,\gamma\beta,
\beta\alpha\gamma-\alpha\gamma\beta \rangle_K.$$
In particular, $\alpha^3=\gamma\beta\in K(\Lambda_2')$ from which it follows that
$$T_1(\Lambda_2'):=\{x\in\Lambda_2'\,\mid\,x^3\in K(\Lambda_2')\} = J(\Lambda_2')$$
where $J(\Lambda_2')$ denotes the ideal generated by the arrows of the quiver. 
(In fact, the basis elements other than $\alpha$ 
in $J(\Lambda_2')$ clearly 
belong to $T_1(\Lambda_2')$ since their third powers vanish.)

For the first K\"ulshammer ideal we therefore get
$$T_1(\Lambda_2')^{\perp} = (J(\Lambda_2'))^{\perp} = \soc(\Lambda_2'),$$
that is, the sequence of K\"ulshammer ideals for $\Lambda_2'$ has length 2.

Secondly, we consider the nonstandard algebra $\Lambda_2$. Now, 
$\beta\gamma=\beta\alpha\gamma$ is nonzero, and the commutator subspace is
$$K(\Lambda_2)=\langle \gamma,\beta,\alpha\gamma,\beta\alpha,\gamma\beta-\beta\alpha\gamma,
\beta\alpha\gamma-\alpha\gamma\beta \rangle_K.$$
In particular, $\alpha\not\in T_1(\Lambda_2):=\{x\in\Lambda_2\,\mid\,x^3\in K(\Lambda_2)\}.$
But then $\alpha^3\in Z(\Lambda_2)$ 
is orthogonal (with respect to the form $\Psi$ above)
to all elements in $T_1(\Lambda_2)$.
Thus, the sequence of K\"ulshammer
ideals takes the form
$$\soc(\Lambda_2) \underbrace{\mbox{~~$\subset $~~}}_{1} T_1(\Lambda_2)^{\perp}
=\langle \alpha^3,\alpha^4,\beta\alpha\gamma \rangle_K
\underbrace{\mbox{~~$\subset $~~}}_{2} Z(\Lambda_2). 
$$
where the numbers under the inclusion signs denote codimensions. 
Since the sequences of K\"ulshammer ideals for the algebras $\Lambda_2'$
and $\Lambda_2$ have different lengths (and different codimensions), we can 
conclude from Theorem \ref{th:zimmermann}
that $\Lambda_2'$ and $\Lambda_2$ are not derived equivalent.  
  
~~~\smallskip
  
\item[{(2)}]
Let $\charact K = 2$, and $\lambda, \mu \in K \setminus \{0,1\}$. 
For both types of algebras $\Lambda_3(\lambda)$ and $\Lambda_3(\mu)'$,
a $K$-linear basis
consisting of the different nonzero paths is given by
$${\mathcal B} = \{e_1,e_2,\alpha,\beta,\sigma,\gamma, \alpha^2, \alpha\sigma,
\beta^2, \gamma\alpha, \alpha^3, \beta^3\}.
$$
In fact, note that for $\Lambda_3(\mu)'$ we have that
$\alpha^4 = (\sigma\gamma)^2 = \mu \sigma\beta^2\gamma
= \mu \alpha\sigma\gamma\alpha
= \mu\alpha^4$, from which $\alpha^4=0 $ follows since $\mu\neq 1$. 
Similarly, one checks that $\gamma\alpha^2=\beta\gamma\alpha=\gamma\sigma\gamma=0$ 
and $\alpha^2\sigma=\alpha\sigma\beta=\sigma\gamma\sigma=0$ in $\Lambda_3(\mu)'$,
and that
$\beta^4=0$ in $\Lambda_3(\lambda)$ and 
$\Lambda_3(\mu)'$.

For both types of algebras $\Lambda=\Lambda_3(\lambda)$ and 
$\Lambda=\Lambda_3'(\mu)$, the center has a basis of the form
$$Z(\Lambda)= \langle 1,\alpha+\beta,\alpha^2,\beta^2,\alpha^3,\beta^3\rangle_K.
$$
It follows from \cite[Section 4]{BiS1} (or can be checked directly) 
that a symmetrizing form 
$\Psi'$ on the standard algebra $\Lambda_3(\mu)'$  
is given by assigning 1 to $\alpha^3$, assigning $\mu^{-1}$ to $\beta^3$,
and 0 to all remaining elements of the basis $\mathcal{B}$.

Note that we can not use the same form for the nonstandard algebra
$\Lambda_3(\lambda)$, because $0=\Psi'(\alpha^2)=\Psi'(\sigma\gamma+\alpha^3)
= 0+1=1$. 
Instead, a symmetrizing form $\Psi$ for $\Lambda_3(\lambda)$ is given by 
assigning 1 to $\alpha^3$ and to $\alpha^2$, assigning $\lambda^{-1}$ to $\beta^3$,
and 0 to all remaining elements of the basis $\mathcal{B}$.

The commutator subspaces for both $\Lambda=\Lambda_3(\lambda)$ and 
$\Lambda=\Lambda_3(\mu)'$
are generated as vector space as follows
$$K(\Lambda) = \langle \sigma, \gamma, \alpha\sigma, \gamma\alpha, 
\gamma\sigma-\sigma\gamma, \gamma\sigma\beta-
\sigma\beta\gamma
\rangle_K.$$
%=\lambda\beta^3-\alpha^3 
Note that the only crucial difference in the relations is that
$\alpha^2=\sigma\gamma$ for the standard algebras $\Lambda_3'(\mu)$,
whereas $\alpha^2=\sigma\gamma+\alpha^3$ for the nonstandard algebras 
$\Lambda_3(\lambda)$. This means that the element $\gamma\sigma-\sigma\gamma$
in the commutator space
equals $\mu\beta^2-\alpha^2$ for $\Lambda_3'(\mu)$, and
$\lambda\beta^2-\alpha^2+\alpha^3$ for $\Lambda_3(\lambda)$, respectively.
This has consequences for the spaces 
$T_1(\Lambda) = \{x\in \Lambda\,\mid\,x^2\in K(\Lambda)\}.$
Namely, let $\sqrt{\mu}$ be the unique square root
in $K$
(recall that $K$ is algebraically closed, and of characteristic 2).
Then we get for the standard algebras
that $\alpha+\sqrt{\mu}\,\beta\in T_1(\Lambda_3'(\mu))$,
whereas there is no analogous element in $ T_1(\Lambda_3(\lambda))$.
More precisely, we have
$$T_1(\Lambda_3'(\mu)) = \langle \alpha+\sqrt{\mu}\,\beta,\sigma,
\gamma, \alpha^2, \alpha\sigma, \beta^2, \gamma\alpha, \alpha^3, \beta^3 
\rangle_K
$$
and 
$$T_1(\Lambda_3(\lambda)) = \langle \sigma,
\gamma, \alpha^2, \alpha\sigma, \beta^2, \gamma\alpha, \alpha^3, \beta^3 
\rangle_K .
$$
From this we can determine the first K\"ulshammer ideals and
their codimensions as follows. For the standard algebras we obtain
(using the symmetrizing form $\Psi'$ above) 
$$\soc(\Lambda_3'(\mu)) \mbox{~~$\underbrace{\subset}_{1}$~~}
T_1(\Lambda_3'(\mu))^{\perp} = \langle \alpha^2+\sqrt{\mu}\,\beta^2,
\alpha^3,\beta^3 \rangle_K \mbox{~~$\underbrace{\subset}_{3}$~~}
Z((\Lambda_3'(\mu)).$$
In fact, $\Psi'((\alpha^2+\sqrt{\mu}\,\beta^2)(\alpha+\sqrt{\mu}\,\beta))
=\Psi'(\alpha^3+\mu\beta^3)=1+\mu\mu^{-1}=0$
since $\charact K=2$; on the other hand, $\alpha^2\not\in T_1(\Lambda_3'(\mu))^{\perp}$
because $\alpha+\sqrt{\mu}\,\beta\in T_1(\Lambda_3'(\mu))'$, similarly for
$\beta^2$. 
\smallskip

For the nonstandard algebras we get
(using the symmetrizing form $\Psi$ above)
$$\soc(\Lambda_3(\lambda)) \mbox{~~$\underbrace{\subset}_{2}$~~}
T_1(\Lambda_3'(\lambda))^{\perp} = \langle \alpha^2,\beta^2,
\alpha^3,\beta^3 \rangle_K \mbox{~~$\underbrace{\subset}_{2}$~~}
Z((\Lambda_3(\lambda)).$$
Because of the different codimensions we can again conclude by Theorem 
\ref{th:zimmermann} that
$\Lambda_3(\lambda)$ and $\Lambda_3'(\mu)$ are not derived equivalent,
as claimed.
 
~~~\smallskip 
 
\item[{(3)}]
Let $\charact K = 2$. For both of the algebras $\Lambda_5$ and $\Lambda_5'$,
a $K$-linear basis consisting of the different nonzero paths is given by
$${\mathcal B}=\{ e_1,e_2,e_3,\alpha,\beta,\gamma,\delta,\sigma,\alpha\gamma,
\alpha^2=\gamma\beta, \beta\alpha,\alpha^3=\delta\sigma,\sigma\delta,
\beta\alpha\gamma\}.
$$
Also for both of the algebras $\Lambda=\Lambda_5$ and $\Lambda=\Lambda_5'$, 
the center is as vector space generated as follows
$$Z(\Lambda)=\langle 1,\alpha^2,\alpha^3, \sigma\delta, \beta\alpha\gamma 
\rangle_K.
$$
It follows from \cite[Section 4]{BiS1} (or can be checked directly) that a 
symmetrizing form 
$\Psi'$ on the standard algebra $\Lambda_5'$ 
is given by assigning 1 to the socle elements $\alpha^3$, $\beta\alpha\gamma$,
and $\sigma\delta$, and assigning 
0 to all other elements of the basis $\mathcal{B}$.

Note that one can not use the same symmetrizing form 
for the nonstandard algebra $\Lambda_5$ since 
$\Psi'(\gamma\beta)=0\neq 1=\Psi'(\beta\alpha\gamma) =\Psi'(\beta\gamma)$,
i.e. the form is not symmetric on $\Lambda_5$. 

In fact, a symmetrizing form $\Psi$ for $\Lambda_5$ is given 
by assigning 1 to 
$\alpha^3$, $\beta\alpha\gamma$, $\sigma\delta$, and
also to $\alpha^2$ and assigning 
0 to all other elements of the basis $\mathcal{B}$.

We shall need the commutator subspaces. It turns out that the only 
crucial difference comes from the fact that $\beta\gamma=0$ in
$\Lambda_5'$, whereas $\beta\gamma=\beta\alpha\gamma$ is nonzero
in $\Lambda_5$. 
%(In fact, note that all other relations are precisely the same.) 
Hence, we get
$$K(\Lambda) = \langle \beta,\gamma,\delta,\sigma,\beta\alpha, \alpha\gamma,
\delta\sigma-\sigma\delta, \beta\alpha\gamma-\alpha\gamma\beta, X
\rangle_K ,
$$
where $X=\gamma\beta$ for $\Lambda_5'$, and $X=\gamma\beta-\beta\gamma$
for $\Lambda_5$, respectively. 
This has the following implication for 
$T_1(\Lambda):=\{x\in\Lambda\,\mid\,x^2\in K(\Lambda)\}$.
Namely,  since
$\alpha^2=\gamma\beta$ we have $\alpha\in T_1(\Lambda_5')$, but
$\alpha\not\in T_1(\Lambda_5)$. More precisely, denoting by
$J(\Lambda)$ the ideal generated by the arrows of the quiver, we obtain
$$T_1(\Lambda_5')= J(\Lambda_5')$$
whereas 
$$T_1(\Lambda_5)\subset J(\Lambda_5)$$
has codimension 1. (In fact, it is easy to see that all other basis elements
of the ideal $J(\Lambda)$ are contained in $T_1(\Lambda)$.)

Now we examine the first K\"ulshammer ideals 
$T_1(\Lambda)^{\perp}$. 
First we consider the standard
algebra $\Lambda_5'$. Here we get
$$T_1(\Lambda_5')^{\perp} =  (J(\Lambda_5'))^{\perp} = \soc(\Lambda_5')
= \langle \alpha^3, \sigma\delta, \beta\alpha\gamma \rangle_K.
$$ 
In fact, $\alpha^2\not\in T_1(\Lambda_5')^{\perp}$ since 
$\alpha\in T_1(\Lambda_5')$ and $\Psi'(\alpha^2\cdot \alpha) = \Psi'(\alpha^3)=1$. 

On the other hand, for the nonstandard algebra $\Lambda_5$ we have 
that $\alpha^2$ is orthogonal to all elements in $T_1(\Lambda_5)$ since
$\alpha\not\in T_1(\Lambda_5)$ (and $\alpha^2\cdot w=0$ for every arrow
$w\neq\alpha$ in the quiver), that is,
$$T_1(\Lambda_5)^{\perp} = \langle \alpha^2,
\alpha^3, \sigma\delta, \beta\alpha\gamma \rangle_K.
$$ 
Since the first K\"ulshammer ideals have different codimensions
inside the center, we can conclude by Theorem \ref{th:zimmermann}
that $\Lambda_5$ and $\Lambda_5'$ are not derived equivalent, as claimed. 

\item[{(4)}]
  Let $\charact K = 2$. The algebras $\Lambda_9$ and 
  $\Lambda'_9$ can not be distinguished by their K\"ulshammer ideals.
  Instead we shall use Hochschild cohomology. It has recently been shown 
by D. Al-Kadi in \cite[Theorems 3.1 and 4.1]{AlK} that
  $$
    \dim_K HH^2(\Lambda_9) = 4
    \mbox{ and }
    \dim_K HH^2(\Lambda_9') = 3 .
  $$
  Therefore, applying Theorem \ref{th:1.5}, we conclude 
  (see \cite[Corollary 4.2]{AlK}) that 
  $\Lambda_9$ and $\Lambda_9'$ are not derived equivalent.
\end{enumerate}

\noindent
This completes the proof of the Main Theorem.



\end{document}